\magnification\magstep 1


\input amssym.tex

\font\titlefont cmb16
\font\authorfont cmr14
\font\headingfont cmr12
\font\twelverm cmr12
\font\twelvei cmmi12

\font\reftext cmr9
\font\refit cmti9
\font\refbf cmbx9

\font\addrtext cmcsc9

\def\Hom{{\rm Hom}}
\def\Ext{{\rm Ext}}
\def\Vect{\mathop{\rm Vect}}
\def\id{{\rm id}}
\def\opp{{\rm opp}}
\def\End{\mathop{\rm End}}
\def\END{\mathop{\cal E\mit nd\.}}
\def\Gr{\mathop{\rm Gr}\nolimits}
\def\Ker{\mathop{\rm Ker}}
\def\Im{\mathop{\rm Im}}
\def\ad{\mathop{\rm ad}\nolimits}

\def\E{{\cal E}}
\def\D{{\cal D}}
\def\C{{\cal C}}
\def\B{{\cal B}}
\def\O{{\cal O}}
\def\U{{\cal U}}

\def\g{{\frak g}}

\def\ot{\otimes}
\def\sub{\subset}
\def\d{\partial}

\let\le\leqslant
\let\ge\geqslant

\def\eps{\varepsilon}
\def\ph{\varphi}

\def\:{\colon}
\def\;{,\>}
\def\.{\mskip.5\thinmuskip}
\def\bu{{ {}^{ {}_{ {}_\bullet } }\mskip-1.3\thinmuskip }}

\def\qed{\quad$\square$}
\def\mathqed{\quad\square}

\emergencystretch=1em

\centerline{\titlefont Nonhomogeneous Quadratic Duality and Curvature}

\vskip 0.85 truecm\relax

\centerline{\authorfont L.~E. Positselski}

\vskip 1 truecm\relax

\centerline{\headingfont Introduction}
\vskip 0.45 truecm\relax

 A {\it quadratic algebra} is a graded algebra with generators
of degree~$1$ and relations of degree~$2$.
 Let $A$ be a quadratic algebra with the space of generators $V$ and
the space of relations $I\sub V\ot V$.
 The classical quadratic duality assigns the quadratic algebra $A^!$
with generators from $V^*$ and the relations $I^\perp\sub V^*\ot V^*$
to the algebra~$A$.
 According to the classical results of Priddy and
L\"ofwall~[1, 3], $\,A^!$ is isomorphic to the subalgebra of
$\Ext^*_A(k,k)$ generated by $\Ext_A^1(k,k)$.
 Priddy called an algebra $A$ a {\it Koszul algebra} if this
subalgebra coincides with the whole of $\Ext_A^*(k,k)$.
 Koszul algebras constitute a wonderful class of quadratic algebras,
which is closed under a large set of operations, contains the main
examples, and perhaps admits a finite classification.

 In this paper, we propose an extension of the quadratic duality to
the nonhomogeneous case.
 Roughly speaking, a {\it nonhomogeneous quadratic algebra} (or 
a {\it quadratic-linear-scalar algebra}, a QLS-{\it algebra})
is an algebra defined by (generators and) nonhomogeneous relations
of degree~$2$.
 A {\it quadratic-linear algebra} (QL-{\it algebra}) is an algebra
defined by nonhomogeneous quadratic relations without the scalar parts;
in other words, it is an augmented QLS-algebra.
 The precise definition takes into account the fact that
a collection of nonhomogeneous relations does not necessarily
``make sense'' (its coefficients must satisfy some equations;
the Jacobi identity is a classical example).

 The dual object for a QL-algebra is~[6] a quadratic
DG-algebra~[7].
 The dual object for a QLS-algebra is a set of data which we call
a {\it quadratic {\rm CDG-}algebra} (``curved''), defined up to
an equivalence.
 The classical Poincar\'e--Birkhoff--Witt theorem on the universal
enveloping algebra structure~[9] finds its natural place in
this context as a particular case of the fact that {\it every
Koszul {\rm CDG-}algebra corresponds to a {\rm QLS-}algebra}.

 Remarkable examples of nonhomogeneous quadratic duality are provided
by Differential Geometry.
 The algebra of differential operators on a manifold may be considered
as a QL-algebra defined by the commutation relations for the vector fields.
 The dual object for this algebra is the de~Rham complex;
the corresponding equivalence of the categories of modules
is constructed in~[4].
 The algebra of differential operators in a vector bundle is
a QLS-algebra.
 The dual object is the algebra of differential forms with coefficients
in linear operators in this bundle and with the exterior differential
defined by means of a connection.
 Its square is not zero---it is equal to the commutator with
the curvature; the equivalence relation mentioned in
the previous paragraph corresponds to changing the connection.
 Thus, the curvature corresponds to the scalar part of the relations.

 A question arises about the obstructions to existence of
a QL-algebra structure (i.e., of ``a flat connection'') on
a given QLS-algebra.
 In the present paper we construct obstructions of this kind
generalizing the Chern classes of vector bundles~[8] and (which is
less evident) the Chern--Weil classes of principal $G$-bundles.
 Our analogues of the secondary characteristic classes~[5]
form the {\it Chern--Simons functor} on the category of CDG-algebras.

 The author is grateful to A.~B.~Astashkevich, R.~V.~Bezrukavnikov,
M.~V.~Finkelberg, V.~A.~Ginzburg, A.~E.~Polishchuk, and V.~S.~Retakh
for numerous and very helpful discussions, and to A.~A.~Kirillov and
A.~N.~Rudakov for their constant attention to~his~work.
{\uchyph=0\par}

\medskip
{\bf Conventions and notation.}
 An {\it algebra} is an associative algebra with a unit~$e$ over
a fixed ground field~$k$; all graded and filtered algebras are
assumed to be locally finite-dimensional.
 The symbol $[a,b]$ denotes the supercommutator
$ab-(-1)^{\tilde{a}\tilde{b}}ba$.

\vskip 0.75 truecm plus 0.1 truecm\relax
\centerline{\headingfont \S1. Definitions}
\vskip 0.4 truecm\relax

{\bf Definition~1.}  A {\it weak {\rm QLS-}algebra} is
an algebra $A$ together with a {\it subspace of generators}
$W\sub A$, $\,e\in W$, satisfying the following conditions.
 Let $V$ be a hyperplane in $W$ complementary to $k\cdot e$,
$\,T(V)$ the tensor algebra, $T_n(V)$ the subspace of elements
of degree~$\le n$, and $J$ the kernel of the natural projection
$T(V)\to A$.
 It is required that \par
  1) an algebra $A$ be generated by its subspace~$W$; \par
  2) the ideal $J$ be generated by its subspace $J_2=J\cap T_2(V)$.
\par\noindent
 (Clearly, this conditon does not depend on the choice of~$V$.)
 In this case, the {\it underlying quadratic algebra} $A^{(0)}$ is
defined by the generators $V\simeq W/k\cdot e$ and the relations
$I=J_2\bmod J_1\sub V\ot V$.

\medskip
{\bf Definition~2.}  A QLS-{\it algebra} is an algebra $A$
together with a filtration $F$: $\,0\sub F_0A\sub F_1A\sub F_2A
\sub\cdots\sub A$, $\,\bigcup F_iA=A$, $\,F_0A=k$,
$F_iA\cdot F_jA\sub F_{i+j}A$, such that the associated graded algebra
$\Gr_FA$ is quadratic.

\smallskip
 A (weak) QLS-algebra $A$ is said to be {\it Koszul\/} if $A^{(0)}$ is
a Koszul algebra.
 The second of our two definitions implies the first one for $W=F_1A$.
 We shall show in~3.3 that in the Koszul case these definitions are
equivalent.

\medskip
{\bf Definition~3.}  A (weak) QL-{\it algebra} is a (weak)
QLS-algebra together with an augmentation (a ring homomorphism)
$\eps\:A\to k$; let $A_+=\Ker\eps$ denote the augmentation ideal.

\medskip
 A {\it morphism of weak {\rm QLS-}algebras} is an algebra
homomorphism preserving the subspace~$W$.
 A {\it morphism of\/ {\rm QLS-}algebras} is an algebra homomorphism
preserving the filtration.
 A {\it morphism of weak {\rm QL-}algebras} is an algebra
homomorphism preserving $W$ and the augmentation.
 This defines the {\it category ${\cal WQLS}$ of weak
{\rm QLS-}algebras} and its full subcategory ${\cal QLS}$,
as well as the {\it category ${\cal WQL}$ of weak
{\rm QL-}algebras} and its full subcategory~${\cal QL}$.
 The {\it categories of Koszul} QLS- and QL-{\it algebras}
are denoted by ${\cal KLS}$ and~${\cal KL}$, respectively.

\medskip
 Now we define the dual objects.

\medskip
{\bf Definition~4.}
 A DG-{\it algebra} $B^\bu$ is a graded algebra (with upper indices)
together with a derivation $d\:B\to B$ of degree~$+1$
$\,$($d(a\cdot b)=d(a)\cdot b+(-1)^{\tilde{a}}a\cdot d(b)$) such that
$d^2=0$.
 A {\it morphism of\/ {\rm DG-}algebras} is a homomorphism of graded
algebras commuting with the derivations.
 We assume below that $B^i=0$ for $i<0$.
 A DG-algebra $B^\bu$ is said to be {\it quadratic} ({\it Koszul\/})
if the algebra $B$ is quadratic (Koszul).
 Let ${\cal QDG}$ and ${\cal KDG}$ be the {\it categories of
quadratic} and {\it Koszul\/ {\rm DG-}algebras}, respectively.

\medskip
{\bf Definition~5.}
 A CDG-{\it algebra} is a triple $\Psi=(B,d,h)$, where $B$ is
a graded algebra, $d$~is a derivation of $B$ of degree~$+1$, and
$h\in B^2$, such that \par
 1) $d^2=[h,{}\cdot{}]$, \par
 2) $d(h)=0$. \par\noindent
In the sequel we assume that $B^i=0$ for $i<0$.
 A CDG-algebra $\Psi$ is called {\it quadratic} ({\it Koszul\/})
if the algebra $B$ is quadratic (Koszul).

\medskip
{\bf Definition~6.}
 A {\it morphism of} CDG-{\it algebras} $\ph\:\Psi\to\Psi'$ is
a pair $\ph=(f,\alpha)$, where $f\:B\to B'$ is a homomorphism of
graded algebras and $\alpha\in B'{}^1$, satisfying the conditions \par
 1) $d'f(x)=f(dx)+[\alpha,f(x)]$, \par
 2) $h'=f(h)+d'\alpha-\alpha^2$. \par\noindent
The {\it composition} $(f,\alpha)\circ(g,\beta)$ is the morphism
$(f\circ g\;\alpha+f(\beta))$.
The {\it identity morphism} is $(\id,0)$.
This defines the {\it category ${\cal CDG}$ of
{\rm CDG-}algebras} and its full subcategories
${\cal QCDG}$ and ${\cal KCDG}$ {\it of quadratic} and
{\it Koszul {\rm CDG-}algebras}, respectively.
Two CDG-algebras $\Psi$ and $\Psi'$ are called {\it equivalent}
if $B=B'$ and there exists a morphism of the form
$(\id,\alpha)\:\Psi \to \Psi'$, in other words, if
$d'=d+[\alpha,{}\cdot{}]$ and $h'=h+d\alpha+\alpha^2$;
in this case we write $\Psi'=\Psi(\alpha)$.

\vskip 0.75 truecm plus 0.1 truecm\relax
\centerline{\headingfont \S2. Duality Functor}
\vskip 0.4 truecm\relax

{\bf 2.1.}  Let $A$ be a weak QLS-algebra with the space
of generators $W\sub A$.
 Set $B=A^{(0)}{}^!$; we will denote by upper indices the grading
on~$B$.
 Choose a hyperplane $V\sub W$ complementary to $k\cdot e$ in~$W$;
we have $J_2\sub k\oplus V\oplus V\ot V$.
 Note that $J_2\cap (k\oplus V)=0$ and $J_2\bmod (k\oplus V)=I$.
 Thus, $J_2$ can be represented as the graph of a linear map
$I\to k\oplus V$, which we will denote by $(-h\;-\ph)$, where
$h\in I^*\simeq B^2$ and $\ph\:I\to V$.
 Let $d_1=\ph^*\:B^1\to B^2$; then the relations in the algebra $A$
can be written in the form
$$
 p+\ph(p)+h(p)=0, \quad p\in I\sub V\ot V, \quad \ph=d_1^*.
 \eqno{(*)}
$$

\medskip
{\bf 2.2. Proposition.}  {\it The map $d_1\:B^1\to B^2$ can be
extended to a derivation~$d$ of the algebra~$B$.
 The triple $(B,d,h)$ is a {\rm CDG}-algebra.}

\medskip
{\bf Proof.} 
 Tensoring the relation~($*$) by $V$ on the left and on
the right, we obtain
$$
 q+\ph_{12}(q)+h_{12}(q) = 0 = q + \ph_{23}(q) + h_{23}(q)
 \pmod J
$$
for any $q\in V\ot I\cap I\ot V$, whence
$-\ph_{12}(q)+\ph_{23}(q)=h_{12}(q)-h_{23}(q)$ $\,({\rm mod}\ J)$.
 The latter equation implies
$$
 \ph_{12}(q)-\ph_{23}(q)\in I  \eqno{(1)}
$$
and
$$
 \ph(-\ph_{12}(q)+\ph_{23}(q)) + h(-\ph_{12}(q)+\ph_{23}(q))
 + h_{12}(q)-h_{23}(q) = 0 \pmod J.
$$
 Since the second summand lies in~$k$, while the first and
the third one belong to $V$, we have
$$
 \eqalignno{
 \ph(\ph_{12}(q)-\ph_{23}(q)) &= h_{12}(q)-h_{23}(q), &(2) \cr
   h(\ph_{12}(q)-\ph_{23}(q)) &= 0.  &(3) \cr}
$$
Dualizing (1), (2), and~(3) and taking into account the fact
that the operator $d_2=(\ph_{12}-\ph_{23})^*$ continues~$d_1$
by the Leibniz rule and that $(h_{12}-h_{23})^*=[h,{}\cdot{}]$,
we obtain, respectively, the equations
($\bmod I^\perp\ot B^1 + B^1\ot I^\perp$):
$$
 d_2(I^\perp)=0, \quad d_2\circ d_1=[h,{}\cdot{}], \quad d_2(h)=0.
$$
 The first equation means that~$d_1$ can be extended to~$B$;
the second and the third one are equivalent to
the CDG-algebra axioms.~\qed

\medskip
{\bf 2.3. Proposition.} {\it Let $V$ and $V'$ be two direct
complements to $k\cdot e$ in~$W$, and let $\Psi$ and $\Psi'$
be the corresponding {\rm CDG-}algebras.
 Then $\Psi'=\Psi(\alpha)$ for a uniquely defined element
$\alpha\in B^1$.}

\medskip
{\bf Proof.}  Suppose that for $\alpha\in B^1=V^*$ we have
$V'=\{v+\alpha(v),\ v\in V\}$, and let $\ph$ and~$h$ 
correspond to the complementary hyperplane~$V$.
 Then for any $p\in I\sub V\ot V$ we have
$$
 \eqalign{
 0 = p+\ph(p)+h(p)=[p&+\alpha_1(p)+\alpha_2(p)+\alpha\ot\alpha(p)] \cr
                     &+[\ph(p)-\alpha_1(p)-\alpha_2(p)+
   \alpha(\ph(p)-\alpha_1(p)-\alpha_2(p))] \cr
                     &+[h(p)-\alpha(\ph(p))+\alpha\ot\alpha(p)]. \cr}
$$
Thus, the operators
$$
 \ph'=\ph-\alpha_1-\alpha_2, \qquad h'=h-\alpha\circ\ph+\alpha\ot\alpha
$$
correspond to the choice of the direct complement $V'$.
 Dualizing and using the fact that
$(\alpha_1+\alpha_2)^*=[\alpha,{}\cdot{}]$ and
$\alpha\circ\ph=d(\alpha)$, we obtain
$d'_1=d_1-[\alpha,{}\cdot{}]$ and $h'=h-d(\alpha)+\alpha^2$.~\qed

\medskip
{\bf 2.4. Proposition.}  {\it The construction of Subsection~{\rm 2.1}
defines a fully faithful contravariant functor
$\D=\D_{\rm QLS}\:{\cal WQLS}\to{\cal QCDG}$.}

\medskip
{\bf Proof.}  To define the functor on objects, choose
the subspace $V$ arbitrarily for every weak QLS-algebra.
 The natural isomorphism
$$
 \{f\in\Hom(k\oplus V\;k\oplus V'):f|_k=\id\}\simeq
 \Hom(V,V')\oplus V^*,
$$
together with Proposition~2.3, allows to define it on morphisms.
 It is obviously fully faithful.~\qed

\medskip
{\bf 2.5.}  Let $A$ be a weak QL-algebra with the augmentation
ideal~$A_+$.
 Set $V=A_+\cap F_1$.
 Then it is easy to see that $h=0$, and we obtain a DG-algebra
$(B,d)$.

\medskip
{\bf Proposition.}  {\it This defines a fully faithful contravariant
functor $\D=\D_{\rm QL}\:{\cal WQL}\to{\cal QDG}$.}

\medskip
{\bf 2.6.}  Conversely, let $\Psi$ be a quadratic CDG-algebra.
Set $V=B_1^*$, $\,I=B_2^*\rightarrowtail V\ot V$,
let $A=A(\Psi)$ be the algebra with generators from $V$ and
relations~($*$), let $W$ be the image of $k\oplus V$ in $A$,
and put $F_nA=W^n$.
 It is easy to see that if $(A(\Psi),W)$ is a weak QLS-algebra
and $A^{(0)}{}^!\simeq B$, then $\D(A,W)=\Psi$.

\medskip
{\bf 2.7. Examples.}
1.~Let $\g$ be a Lie algebra.
 Then the enveloping algebra $U\g$ is a QL-algebra, and any
QL-algebra $A$ for which $\Gr_FA$ is a symmetric algebra can be
obtained in this way.
 The dual DG-algebra $(\Lambda\!{}^\bu\g^*,d)$ is the standard 
cohomological complex of the Lie algebra~$\g$.
 The QLS-algebras $A$ for which $\Gr_FA$ is a symmetric algebra
correspond to central extensions of Lie algebras:
the algebra $A=U\widetilde{\g}/(1-c)$ is assigned to a central
extension $0\to k\cdot c\to\widetilde{\g}\to\g\to0$.
 The dual object is $\Psi=(\Lambda\!{}^\bu\g^*,d,h)$, where $h$
is the cocycle of the central extension.

2.~The Clifford algebra $\{vw+wv=Q(v,w),\ v,\,w\in V\}$ is
a QLS-algebra, and all QLS-algebras with $\Gr_FA=\Lambda\!{}^\bu V$
are Clifford algebras.
 The dual object is $\Psi=(S^\bu V,0,Q)$.
 All QL-algebras with $\Gr_FA=\Lambda\!{}^\bu V$ have the form
$\{vw+wv=\lambda(v)w+\lambda(w)v\}$, $\,\lambda\in V^*$.

3. Let $A=k\oplus A_+$ be a (finite-dimensional) augmented algebra.
 Let us endow $A$ with the structure of a QL-algebra by
setting $F_iA=A$ for $i\ge1$.
 Then the dual DG-algebra is the reduced cobar-construction
for~$A$, $\,\D(A,F)=\C^\bu(A)=\sum^\oplus A_+^{*\ot n}$.

\medskip
{\bf 2.8. Remark.}  Under the quadratic duality, commutative algebras
correspond to universal enveloping algebras of Lie algebras.
 In particular, we have \par
a)~a duality between commutative QLS-algebras and quadratic Lie
CDG-algebras, and \par
b)~a duality between Lie QL-algebras and supercommutative
quadratic DG-algebras.

\vskip 0.75 truecm plus 0.1 truecm\relax
\centerline{\headingfont \S3. Bar Construction}
\vskip 0.4 truecm\relax

{\bf 3.1. Bar-complex for CDG-algebras.}  The following construction
is due to A.~E.~Polishchuk.
 Let $\Psi=(B,d,h)$ be a CDG-algebra, $B_i=0$ for $i<0$, $\,B_0=k$.
 Put $\B(\Psi)=\sum^{\oplus\,\infty}_{n=0}B_+^{\ot n}$, where
$B_+=\sum_{i=1}^\infty B_i$, and, denoting by
$(b_1\mid b_2\mid\ldots\mid b_n)$ the element
$b_1\ot b_2\ot\cdots\ot b_n\in\B(\Psi)$, endow $\B(\Psi)$ with
a coalgebra structure:
$$
 \Delta(b_1\mid b_2\mid\ldots\mid b_n) =
 \sum_{k=0}^n(b_1\mid\ldots\mid b_k)\ot (b_{k+1}\mid\ldots\mid b_n).
$$
 There are two gradings on $\B(\Psi)$, namely, the internal and
the homological ones:
$$
 \deg_i(b_{i_1}\mid b_{i_2}\mid\ldots\mid b_{i_n}) = i_1+i_2+\cdots+i_n,
 \qquad \deg_h(b_{i_1}\mid b_{i_2}\mid\ldots\mid b_{i_n}) = n
$$
for $b_{i_j}\in B_{i_j}$; set $\B^k=\{b: \deg_ib-\deg_hb=k\}$.
 Let us define the differentials~$\d$, $d$, and~$\delta$ on $\B^\bu$
(of bidegrees $(0\;-1)$, $(1,0)$, and $(2,1)$, respectively)
by the formulas
$$
 \eqalign{
 \d(b_{i_1}\mid\ldots\mid b_{i_n}) &=
 \sum_{k=1}^{n-1} (-1)^{i_1+\ldots+i_k+k-1}
 (b_{i_1}\mid\ldots\mid b_{i_k}b_{i_{k+1}}\mid\ldots\mid b_{i_n}), \cr
 d(b_{i_1}\mid\ldots\mid b_{i_n}) &=
 \sum_{k=1}^n (-1)^{i_1+\ldots+i_{k-1}+k-1}
 (b_{i_1}\mid\ldots\mid d(b_{i_k})\mid\ldots\mid b_{i_n}), \cr
 \delta(b_{i_1}\mid\ldots\mid b_{i_n}) &=
 \sum_{k=1}^{n+1} (-1)^{i_1+\ldots+i_{k-1}+k-1}
 (b_{i_1}\mid\ldots\mid b_{i_{k-1}}\mid h\mid b_{i_k}
 \mid\ldots\mid b_{i_n}). \cr}
$$
 It is straightforward to check that~$\d$, $d$, and~$\delta$ are
superderivations of the coalgebra $\B$ and $(d+\d+\delta)^2=0$.
 Let $(\C_\bu(\Psi),D)$ be the dual DG-algebra to
the DG-coalgebra $\B^\bu(\Psi)$:
$$
 \C_\bu(\Psi)=\sum_{n=1}^\infty\!{}^{{}^{\scriptstyle\oplus}}
 \bigg(\sum_{i=1}^\infty\!{}^{{}^{\scriptstyle\oplus}}B_i^*\bigg)
 ^{\ot n}, \qquad D=(\d+d+\delta)^*.
$$

\medskip
{\bf 3.2. Definition.}
 The {\it bar-cohomology algebra\/} of a CDG-algebra $\Psi$ is 
the homology algebra of the DG-algebra $\C_\bu(\Psi)$, \
$H_\bu^b(\Psi)=H_\bu(\C(\Psi),D)$.

\medskip
{\bf Proposition} (L\"ofwall's subalgebra theorem).
{\it If\/ $\Psi$ is quadratic, then $H^b_0(\Psi)$ is isomorphic to
the algebra $A(\Psi)$ constructed in\/~{\rm 2.6.}
 The filtration $F_n=F_1^n$ on $A$ is induced by the\/
$\deg_i$-filtration of the cobar-complex\/ $\C_\bu(\Psi)$.}

\medskip
 The proof is immediate. \qed

\medskip

{\bf Corollary.}
{\it If $A$ is a weak\/ {\rm QLS-}algebra, then $H^b_0(\D(A))=A$.}

\medskip
{\bf 3.3.  Poincar\'e--Birkhoff--Witt theorem.}

\medskip
{\bf Theorem.}
{\it Let\/ $\Psi=(B,d,h)$ be a Koszul {\rm CDG-}algebra.
Then there is an isomorphism\/ $\Gr A(\Psi)\simeq B^!$.}

\medskip
{\bf Proof.} (We shall see that the Koszul condition can be weakened
to the requirement that $\Ext^i_B(k,k)_{i+1}=0$ for all~$i$, or,
equivalently $\Ext^3_{B^!}(k,k)=0$ for $i\ge4$.)
 There is a spectral sequence $E^1_{p,q}=\Ext_B^{-q}(k,k)_p
\Rightarrow H^b_{p+q}(\Psi)$ induced by the $\deg_i$-filtration
on $\C_\bu(\Psi)$.
 Since $B$ is Koszul, we have $E^1_{p,q}=0$ for $p+q\ne0$, and
$E^r_{p,q}$ degenerates at the term~$E^1$.
 Therefore, $\Gr H^b_\bu(B,d,h)=\Ext_B(k,k)=B^!$. \qed

\medskip
{\bf Corollary.}
{\it Any Koszul weak\/ {\rm QLS-(QL-)}algebra $A$ is a\/
{\rm QLS-(QL-)}algebra, and $H^b(\D(A))=A$.
The duality functors
$$
 \D_{\rm QL}\:{\cal KL}\to {\cal KDG} \quad \hbox{and}
 \quad \D_{\rm QLS}\:{\cal KLS}\to {\cal KCDG}
$$
are antiequivalences of categories.}

\medskip
{\bf Proof.}
 Apply Subsection~2.6, the proof of the theorem, and the fact
that the algebras $B$ and $B^!$ are Koszul simultaneously. \qed

\medskip
{\bf 3.4.} Without the Koszul condition the statement of
Theorem~3.3 fails.
 A counterexample~[6]: the relations
$$
 xy=x+y, \quad x^2+yz=z
$$
imply $yz=zy$, although they have the form~($*$) for a certain
DG-algebra.

\vskip 0.75 truecm plus 0.1 truecm\relax
\centerline{\headingfont\textfont0=\twelverm\textfont1=\twelvei
\S4. An Example: $D$-$\Omega$-Duality}
\vskip 0.4 truecm\relax

 Strictly speaking, these examples do not keep within our scheme,
and we shall only show that they are similar to it (however,
the scheme can be extended to include them).

\medskip
{\bf 4.1.}  Let $M$ be a smooth manifold, $\O(M)$ the ring of
smooth functions on~$M$, $\,\E$ a vector bundle on~$M$,
$\,D(M,\E)$ the ring of differential operators in~$\E$, and
$F_nD(M,\E)$ the subspace of operators of degree at most~$n$.
 The equation
$$
 \Gr_F D(M,\E)=\End\E\ot_{\O(M)}S^*_{\O(M)}(\Vect(M))
$$
allows us to consider $\Gr_FD(M,\E)$ as a ``quadratic algebra over
$\End\E$'' and $D(M,\E)$ as a QLS-algebra; then
$\Gr_FD(M,\E)^!=\Omega^*(M\;\END\E)$ is the algebra of differential
forms on $M$ with coefficients in $\END\E$.

 In order to construct a direct complement $V$ to $\End\E$ in
$F_1D(M,\E)$, we choose a connection $\nabla$ on $\E$, define
an embedding
$$
 i\:\End\E\ot\Vect M\to F_1D(M,\E),
 \qquad i(a\ot v)=a\circ \nabla_v,
$$
and put $V=\Im i$.
 It is easy to see that all left $\End\E$-invariant direct complements
to $\End\E$ in $F_1D(M,\E)$ can be obtained in this way.

 Let $d^\nabla$ be the de~Rham differential on $\Omega^*(M\;\END\E)$
defined by means of the connection on $\END\E$ induced by~$\nabla$,
and let $h^\nabla\in\Omega^*(M\;\END\E)$ be the curvature of
the connection~$\nabla$.
 Comparing the relation
$$
 \nabla_X\nabla_Y - \nabla_Y\nabla_X = 
 \nabla_{[X,Y]} + h^\nabla(X,Y)
$$
with the formula~($*$) and taking into account the relationship
between $[{}\cdot{},{}\cdot{}]$ and $d_1$, we conclude that
$\D(D(M,\E)^\opp\;F) = (\Omega^*(M\;\END\E)\;d^\nabla,h^\nabla)$.

\medskip
{\bf 4.2. Principal bundles.}

\medskip
{\bf 4.2.1.}  Let $M$ be a manifold, $G$ a Lie group, $P$ a (right)
principal $G$-bundle over $M$, and $\pi\:P\to M$ the corresponding
projection.
 Let $\g_p$ be the bundle of Lie algebras over $M$ associated with
$P$ by means of the adjoint representation of $G$ (in other words,
the sections of~$\g_p$ are $G$-equivariant vector fields on~$P$),
and let $\U_p$ be the corresponding bundle of enveloping algebras.

\medskip
{\bf Definition.}  The {\it ring of differential operators on
a principal $G$-bundle} is the ring of $G$-equivariant differential
operators on its total space, $D(M,P)=D(P)^G$.
 The filtration~$F$ ``by the order along the base'' on
$D(M,P)$ is defined as follows:
$$
 F_nD(M,P) = \{D\in D(M,P): \ad^{n+1}(\pi^*f)(D)=0
 \ \forall f\in\O(M)\}.
$$

\medskip
{\bf 4.2.2. Proposition.}
$\Gr_FD(M,P)\simeq\U_p\ot_{\O(M)}S^*_{\O(M)}(\Vect(M))$.

\medskip
{\bf Proof.}
 First one has to show that $F_0D(M,P)\simeq\Gamma(\U_p)$.
 Then the isomorphism is defined using the highest symbol operator
$$
 \sigma_n\: F_nD(M,P)\to \Gamma(M\;\U_p\ot S^nT(M)),
$$
$\sigma_n(D)(\xi)=\ad^n(\pi^*f)(D)_m\in\U_{p,\.m}$ for
$m\in M$, $\,\xi\in T^*_m(M)$, $\,f\in\O(M)$, and $d_mf=\xi$.~\qed

\medskip
{\bf 4.2.3.} Let us choose a connection $\nabla$ on the principal
$G$-bundle $P$ and construct a direct complement $V$ to
$F_0D(M,P)$ in $F_1D(M,P)$ as follows:
$V=\langle u\cdot H_\nabla(v), \ u\in\Gamma(\U_p), \
v\in\Vect(M)\rangle$, where $H_\nabla(v)$ is the horizontal
(with respect to~$\nabla$) lifting of the vector~$v$ to~$P$.

 Then $\D_{\rm QLS}(D(M,P)\;F) = (\Omega^*(M,\U_p)\;d^\nabla\;h^\nabla)$,
where $d^\nabla$ is defined by means of the connection $\nabla^u$
on $\U_p$ associated with~$\nabla$, and $h^\nabla\in
\Omega^2(M,\g_p)\subset\Omega^2(M,\U_p)$ is the curvature of~$\nabla$.

\vskip 0.75 truecm plus 0.1 truecm\relax
\centerline{\headingfont \S5. Characteristic Classes}
\vskip 0.4 truecm\relax

 In this section we suppose that the characteristic of the ground
field is equal to~$0$.

\medskip
{\bf 5.1.} Let $\Psi_0=(B,\delta_0,h_0)$ be a CDG-algebra, $[B,B]$ be
the linear subspace generated by the supercommutators in~$B$,
$\,C=B/[B,B]$, and $T\:B\to C$ be the projection.
 It is clear from the Leibniz identity that the operator~$\delta_C$
on $C$ induced by $\delta_0$ is well-defined.
 Notice that $\delta_C^2=0$ and $\delta_C$~does not change when
a CDG-algebra $\Psi$ is replaced by an equivalent one.
 We put $h(\alpha) = h_0 + \delta_0\alpha + \alpha^2$ and
$\delta(\alpha) = \delta_0 + [\alpha,{}\cdot{}]$ for any $\alpha\in B^1$.

\medskip
{\bf 5.2. Main lemma.}
{\it There exist naturally defined differential forms\/ $\omega_n^{(i)}$
on the vector space $B^1$, $\,n=1$, $2$,~\dots, $\,i=0$,~\dots, $n$,
satisfying the following conditions:

\smallskip
\itemitem{\rm(i)} $\omega_n^{(i)}$ is a differential $i$-form with
values in $B^{2n-i}$;
\itemitem{\rm(ii)} $\omega_n^{(0)}=h^n$;
\itemitem{\rm(iii)} $d\omega_n^{(i)}=\delta\omega_n^{(i+1)}$,
$\,i=0$,~\dots, $n-1$, where $d$ is the de~Rham differential;
\itemitem{\rm(iv)} $d\omega_n^{(n)}=0$.
\smallskip}

\medskip
{\bf Proof.}  We put $\omega_n^{(i)}=\sum (d\alpha)^i h^{n-i}$, where
$d\alpha$ is the tautological $1$-form on $B^1$ with values in $B^1$
and $\sum$ denotes the summation over all rearrangements of
the factors.
 Verification is based on the identities $\delta(h)=0$,
$\,dh=\delta\,d\alpha$.  \qed

\medskip
{\bf 5.3. Chern classes.}  Set $c_n=T(h_0^n)\in C^{2n}$.
The {\it Chern classes\/} are the cohomology classes of 
the elements~$c_n$.

\medskip
{\bf Theorem.} a)~$\delta_C(c_n)=0$. \par
b)~The cohomology class of~$c_n$ does not change when
a CDG-algebra $\Psi$ is replaced by an equivalent one.

\medskip
{\bf Proof.} a)~Moreover, $\delta_0(h_0^n)=0$. \par
b)~follows from the equation $d(h^n)=\delta(\omega_n^{(1)})$.  \qed

\medskip
In characteristic~$p$ the theorem remains true for $2n<p$.

\medskip
{\bf 5.4.}  If $\Psi_0=(\Omega^*(M,\END\E)\;d^\nabla\;h^\nabla)$,
then $C^\bu=\Omega^\bu(M)$, the map $T$ is the (matrix) trace,
and we obtain the usual Chern classes.

\medskip
{\bf 5.5. Chern classes: the case of a principal bundle.}

\medskip
{\bf 5.5.1. Lemma.}  Let $A$ be a graded algebra generated by
its graded vector subspace~$W$.
 Then $[A,A]=[W,A]$.

\medskip
{\bf Proof.}  Proceed by induction using the identity
$$
 [ab,c] = [a,bc] + (-1)^{(\tilde b+\tilde c)\tilde a}[b,ca]. \mathqed
$$

\medskip
{\bf Proposition.}  Let $\g$ be a Lie algebra. Then
$$
 (U\g/[U\g,U\g])\simeq (S^\bu\g)_\g.
$$

\medskip
{\bf Proof.}  The map $S^\bu\g\to U\g$, $\,x^n\mapsto x^n$, is
an isomorphism of $\g$-modules, and $U\g/[U\g,U\g]\simeq (U\g)_\g$
by the lemma. \qed

\medskip
{\bf 5.5.2. Proposition.} {\it
{\rm a)}~$\Omega^*(M,\U_p)/[\Omega^*(M,\U_p),\Omega^*(M,\U_p)]\simeq
\Omega^*(M,\U_p/[\U_p,\U_p])$. \par
{\rm b)}~The vector bundle $\U_p/[\U_p,\U_p]$ is trivial with the fiber
$(S^\bu\g)_\g$. \par
{\rm c)}~The form~$c_n$ lies in the space\/ $\Omega^{2n}(M,(S^n\g)_\g)$.
Let $P$ be an invariant polynomial on~$\g$ of degree~$n$.
Then $\langle P,c_n\rangle$ is equal to the Chern--Weil
characteristic form $P(h^\nabla)$ corresponding to $P$ {\rm[8]},
where angle brackets $\langle\,\ ,\ \rangle$ denote the pairing}
$$
 (S^\bu\g^*)^\g\times(S^\bu\g)_\g\to \Bbb R.
$$

\medskip
{\bf Proof.}
a)~If $A$ is a supercommutative algebra, then
$$
 A\ot B/[A\ot B,A\ot B] = A\ot B/[B,B].
$$

b) The adjoint action of $G$ is trivial in $\U_p/[\U_p,\U_p]$.

c) follows from the definition of the isomorphism in
Proposition~5.5.1.

\medskip
{\bf 5.6. The Chern--Simons functor.}

\medskip
{\bf Definition.}  The {\it category\/ $\C2$ of two-term complexes\/}
is defined as follows.
 Its {\it objects\/} are pairs $(C;c)$, where $C=(\delta\:C_1\to C_0)$
is a morphism of vector spaces and $c\in C_0$.
 {\it Morphisms\/} from $(C;c)$ to $(C',c')$ are pairs $(f;c_1')$,
where $f=(f_0,f_1)$, $\,f_i\:C_i\to C_i'$ is a pair of morphisms forming
a commutative square with $\delta$ and~$\delta'$, and $c_1'\in C_1'$
is an element for which $c'-f_0(c) = \delta'c'_1$.
 The {\it composition of morphisms} is defined by the formula
$(f,c_1'')\circ(g,c_1')=(f\circ g\;c_1''+f(c_1'))$.

\medskip
{\bf Construction.}
 The {\it Chern--Simons functor\/} $CS_n\:{\cal CDG}\to \C2$ is
constructed as follows.
 On objects:
$$
 CS_n(\Psi) = (\delta_C\: C^{2n-1}/\delta_C C^{2n-2}\to
 C^{2n}\cap\Ker\delta_C;\mskip\thickmuskip c_n),
$$
where $c_n=T(h_0^n)$.  On morphisms:
$$
 CS_n(f,\alpha) = (f_*,c_n^{(1)}),
 \qquad c_n^{(1)}=c_n^{(1)}(f,\alpha)=\int_\gamma\omega_n^{(1)},
$$
where $f\:\Psi\to\Psi'$, $\,\omega_n^{(1)}$ is the $1$-form
corresponding to the algebra $\Psi'(-\alpha)$, and
$\gamma$~is a smooth path in $B'{}^1$ joining the points $0$
and~$\alpha$.

\medskip
{\bf 5.7. Chern--Simons classes.}
 Let $E=(E,\d)$ be a DG-algebra, and let $\phi\:\Psi\to E$ be
a morphism of CDG-algebras.
 Then $c_n^{(1)}(\phi)\in E^{2n-1}/([E,E]+\d E)$, $\,\d c_n^{(1)}
= -f_*(c_n)$, and when the algebra $\Psi$ is replaced by
an equivalent one the chain $c_n^{(1)}$ changes by an element
from $f(B^{2n-1})$.
 Thus, the class
$$
 c_n^{(1)}\in E^{2n-1}/(f(B)+[E,E]+\d E)
$$
is an invariant of the morphism~$\phi$.

 Let $P$ be a principal $G$-bundle over~$M$.
 When $\Psi=(\Omega^*(M,\U_p)\;d^\nabla\;h^\nabla)$,
$\,(E,\d)=(\Omega^*(P,U\g)\;d)$, $\,f=\pi^*$, and
$\alpha$~is the connection form, one obtains the usual
Chern--Simons classes.

\bigskip

\noindent {\it Corrections made twenty years later:}

1. The Chern classes of 5.4 are more precisely described as
the components of the Chern character (up to factorial factors).
 In other words, they correspond to power sums of symmetric variables
rather than to the elementary symmetric polynomials.

2. Proposition 5.5.2(b) only holds as stated when the Lie group $G$
is connected.
One has to make a separated consideration of invaritant polynomials
for nonconnected Lie groups in this case.

\bigskip

\noindent {\it Notes added twenty years later:}

1. The ``extended scheme including the examples'' promised in \S4 was
indeed worked out (even if not in the most detailed or easily
accessible form) in the auxiliary material to the author's monograph
``Homological algebra of semimodules: Semi-infinite homological algebra
of associative algebraic structures'', Sections 0.4.3--0.4.4 and
11.5--11.6.

\medskip

2. The most important aspect of the CDG-ring theory that was overlooked
in the original 1992-93 paper is that CDG-rings actually form
a {\it $2$-category\/} rather than just a 1-category.
 While CDG-rings themselves describe the curvatures and their 1-morphisms
are responsible for changing connections, the 2-morphisms correpond to
the {\it gauge transformations}.

 Let $\Psi=(B,d,h)$ and $\Psi'=(B',d',h')$ be two CDG-algebras, and let
$(f,\alpha)\:\Psi\to\Psi'$ and $(g,\beta)\:\Psi\to\Psi'$ be two CDG-morphisms
between them.
 A {\it $2$-morphism} $(f,\alpha)\to(g,\beta)$ is an invertible element
of degree zero $z\in B'{}^0$ satisfying the conditions \par
 1) $g(x) = zf(x)z^{-1}$ for all $x\in B$, \par
 2) $\beta = z\alpha z^{-1} + d'(z)z^{-1}$. \par\noindent
 If a pair $(f,\alpha)$ is a morphism of CDG-algebras and $z$ is an invertible
element in $B'{}^0$, then the pair $(g,\beta)$ defined by the above formulas
is also a morphism of CDG-algebras.
 Notice the difference between DG- and CDG-morphisms: while invertible cocycles
of degree zero act by adjunctions on DG-morphisms, invertible cochains of
degree zero act by adjunction on CDG-morphisms.
 
 The 2-category structure on CDG-rings may be possibly used to defined
(quasi-coherent) stacks of CDG-algebras, extending the definitions of
quasi-coherent sheaves of CDG-algebras given in Appendix~B.1 to the author's
memoir ``Two kinds of derived cagegories, Koszul duality, and
comodule-contramodule correspondence'' and Section~1.2 to the preprint
``Cohenent analogoues of matrix factorizations and relative singularity
categories''.
 Under the $D$-$\Omega$ duality, these would correspond to a certain kind of
twisted differential operators (e.g., in the \'etale or analytic topology).

 Finding a 2-category version of the Chern--Simons functor construction
of Subsection~5.6 would be also interesting.

\vskip 0.9 truecm plus 0.1 truecm\relax
\centerline{\headingfont References}
\vskip 0.3 truecm\relax

{\reftext  \baselineskip=0.9\baselineskip
\lineskip=0.9\lineskip  \lineskiplimit=0.9\lineskiplimit

\noindent\phantom{1}[1]
 S.~Priddy.
   {\refit Koszul resolutions.}
Trans.\ Amer.\ Math.\ Soc.\ {\refbf 152}, \#1, p.~39--60, 1970.

\noindent\phantom{1}[2]
 A.~A.~Beilinson, V.~A.~Ginzburg, V.~V.~Schechtman.
   {\refit Koszul duality.}
Journ.\ Geometry and Physics {\refbf 5}, \#3, 317--350, 1988.

\noindent\phantom{1}[3]
 C.~L\"ofwall.
   {\refit On the subalgebra generated by one-dimensional elements
in the Yoneda Ext-algebra.}
Lecture Notes in Math.\ {\refbf 1183}, p.~291--338, 1986.

\noindent\phantom{1}[4]
 M.~M.~Kapranov.
   {\refit On DG-modules over the de Rham complex and the vanishing
cycles functor.}
Lecture Notes in Math.\ {\refbf 1479}, p.~57--86, 1991.

\noindent\phantom{1}[5]
 S.-S.~Chern, J.~Simons.
   {\refit Characteristic forms and geometric invariants.}
Annals of Math.~(2) {\refbf 99}, \#1, p.~48--69, 1974.

\noindent\phantom{1}[6]
 A.~E.~Polishchuk, L.~E.~Positselski.
   {\refit Quadratic algebras}, to appear.

\noindent\phantom{1}[7]
 S.~MacLane.
   {\refit Homology}.
Springer-Verlag, Berlin--New York, 1963.

\noindent\phantom{1}[8]
 J.~W.~Milnor, J.~D.~Stasheff.
   {\refit Characteristic Classes.}
Annals of Math.\ Studies, 76, Princeton University Press,
University of Tokyo Press, 1974.

\noindent\phantom{1}[9]
 J.-P.~Serre.
   {\refit Lie Algebras and Lie Groups.}
Benjamin, New York--Amsterdam, 1965.

\noindent[10]
 F.~A.~Beresin, V.~S.~Retakh.
   {\refit A method of computing characteristic classes of
vector bundles.}
Reports on Math.\ Physics {\refbf 18}, \#3, p.~363--378, 1980.
}

\vskip 0.7 truecm plus 0.1 truecm\relax

\noindent {\addrtext Moscow State University \hfill 1992--93}

\bye